\documentclass[11pt]{amsart}
\usepackage{amscd,amsmath,amsthm,amssymb}
\DeclareMathOperator{\Pic}{Pic}
\DeclareMathOperator{\Aut}{Aut}
\DeclareMathOperator{\Def}{Def}
\DeclareMathOperator{\Image}{Im}

\begin{document}

\newtheorem*{theo}{Theorem A}
\newtheorem{theorem}{Theorem}[section]
\newtheorem{lemma}[theorem]{Lemma}
\newtheorem{proposition}[theorem]{Proposition}
\newtheorem{corollary}[theorem]{Corollary}

\theoremstyle{definition}
\newtheorem{definition}[theorem]{Definition}
\newtheorem{example}[theorem]{Example}

\theoremstyle{remark}
\newtheorem{remark}[theorem]{Remark}


\newcommand{\contr}{{\mspace{1mu}\lrcorner\mspace{1.5mu}}}
\newcommand{\de}{\partial}
\newcommand{\debar}{{\overline{\partial}}}
\newcommand{\desude}[2]{{\dfrac{\de #1}{\de #2}}}

\newcommand{\Morcat}[1]{{\operatorname{Mor}_{\cat{#1}}}}

\newcommand{\ope}[1]{\operatorname{#1}}

\newcommand{\mapor}[1]{{\stackrel{#1}{\longrightarrow}}}
\newcommand{\mapver}[1]{\Big\downarrow\vcenter{\rlap{$\scriptstyle#1$}}}
\newcommand{\implica}[2]{{$[#1\Rightarrow#2]$}}
\newcommand{\biimplica}[2]{{$[#1\Leftrightarrow#2]$}}

\renewcommand{\bar}{\overline}

\newcommand{\Oh}{\mathcal{O}}
\newcommand{\sA}{\mathcal{A}}
\newcommand{\sD}{\mathcal{D}}
\newcommand{\sI}{\mathcal{I}}
\newcommand{\sL}{\mathcal{L}}
\newcommand{\sX}{\mathcal{X}}
\newcommand{\sY}{\mathcal{Y}}

\newcommand{\Cqbar}{\overline{\mathbb{C}^q}^\vee}


\title[Deformations of products and  obstructed
irregular surfaces]{Deformations of products and new examples of obstructed
irregular surfaces}
\author[M. Manetti]{Marco Manetti}

\address{\newline Dipartimento di Matematica ``Guido Castelnuovo'',\hfill\newline
Universit\`a di Roma ``La Sapienza'',\hfill\newline
P.le Aldo Moro 5,\hfill\newline
I-00185 Roma\\ Italy.}
\email{manetti@mat.uniroma1.it}
\urladdr{www.mat.uniroma1.it/people/manetti/}

\date{january 9, 2006}

\begin{abstract}
We determine the base space of the Kuranishi family of some
complete intersection in the product of an abelian variety and a
projective space. As a consequence we obtain new examples of
obstructed irregular surfaces with ample canonical bundle and
maximal Albanese dimension.

Mathematics Subject Classification (2000): 13D10, 14D15.
\end{abstract}

\maketitle

\section*{Introduction}
Define an equivalence relation on analytic singularities generated
by: \emph{if $(X, p)\to(Y, q)$ is a smooth morphism, then $(X,
p)\sim(Y, q)$.} We call the equivalence classes \emph{singularity types}.
It is easy to prove that every analytic singularity is
determined, up to isomorphism,
by its singularity type and the dimension of its Zariski tangent space.\\
For every complex manifold $X$ we denote by  $\Def(X)$ the
deformation space of $X$, i.e. the base space of its semiuniversal
deformation. Thus $\Def(X)$ is determined by its singularity type and by
the dimension of the cohomology group $H^1(X,T_X)$.\\
In the paper \cite{vakil}, Ravi Vakil shows that the moduli space
of regular surfaces satisfies the "Murphy's law". More precisely
he proves that every singularity type defined over $\mathbb{Z}$ is
obtained as deformation space of a regular surface. The methods
used there fail when applied to varieties $X$ with
$H^1(\Oh_X)\not=0$; on the other side there exist in literature
several examples of obstructed surfaces either admitting
irrational pencils (\cite{kas1}, \cite{konno}) or containing nodal
curves
(\cite{BW}, \cite{cat}, \cite{kas2}).\\
The aim of this paper is to extend some of the technical tools
used in \cite{vakil} to irregular surfaces; as a by-product 
we obtain some new easy examples  of obstructed irregular
surfaces.\\

\textbf{Example.} Let $A$ be an abelian surface and $S$ a smooth
surface of general type contained in $A\times \mathbb{P}^1$. Then
$\Def(S)$  has the same singularity
type of the affine cone over
the Segre variety $\mathbb{P}^2\times \mathbb{P}^1\subset \mathbb{P}^5$
(Example \ref{exa.easyexample}).\\

The same ideas can be used to understand the deformation type of
certain complete intersections in
$A\times\mathbb{P}^{n_1}\times\cdots\times \mathbb{P}^{n_k}$. For
example we prove (Theorem \ref{thm.main}):

\begin{theo} Let $A$ be an abelian variety of dimension $q$,  let $D$ be a
sufficiently ample divisor on $A\times\mathbb{P}^{n-1}$
and $S\subset A\times\mathbb{P}^{n-1}$ the intersection of
$m$ generic hypersurfaces homologically equivalent to $D$.\\
If $0\le m\le q+n-3$, then $\Def(S)$ has the same singularity type of the commuting
variety
\[ C(q,\mathfrak{sl}(n,\mathbb{C}))=\{(A_1,\ldots,A_q)\in \mathfrak{sl}(n,\mathbb{C})^{\oplus q}\mid
A_iA_j=A_jA_i\text{\rm\  for every } i,j\}.\]
\end{theo}

Since two matrices in $\mathfrak{sl}(2,\mathbb{C})$ commute if and only if they are linearly
dependent, the commuting variety $C(q,\mathfrak{sl}(2,\mathbb{C}))$ is isomorphic to the affine cone
over the Segre variety $\mathbb{P}^{q-1}\times \mathbb{P}^2\subset \mathbb{P}^{3q-1}$.
By a classical result (\cite{Ger}, \cite{Motz}) the commuting
variety $C(2,\mathfrak{sl}(n,\mathbb{C}))$ is irreducible for every $n$, while
$C(q,\mathfrak{sl}(n,\mathbb{C}))$ is reducible for every pair of positive integers $(q,n)$
such that the product $(q-3)(n-3)$  is sufficiently large (Lemma
\ref{lem.irred}).

The proof of Theorem A is divided in two parts. 
In the first we study deformations of products of K\"{a}hler manifolds (hence the name of this
paper) using the theory of differential graded Lie algebras. 
The main results is the Formality
Theorem \ref{thm.prodotto} which implies in particular that, under suitable cohomological 
condition  on $X,Y$, the deformation space of $X\times Y$ 
is completely determined by the primary obstruction
map of Kodaira and Spencer \cite{KS2}.\\
The second part is essentially a  works about stability and costability theorems; 
more precisely we prove  analogs of Horikawa's theorems \cite{Horidhm1,Horidhm2,Horidhm3}
in some cases where Horikawa's hypothesis are not completely satisfied.\\

The author thanks Corrado De Concini for advices about commuting varieties 
and Edoardo Sernesi 
for useful comments on the first version of this paper.

\section{General notation and basic facts}

We work over the field $\mathbb{C}$ of complex numbers; every 
complex manifold is assumed compact and connected.\\
For every  complex manifold $X$ we denote by:
\begin{itemize}

\item $B_{X}=\oplus_{i} H^{i}(X,\Oh_{X})$ the graded algebra of
the cohomology of the structure sheaf of $X$, endowed with the cup
product $\wedge$.

\item For every holomorphic vector bundle $E$ on $X$ let
$\sA_X^{p,q}(E)$ be  the sheaf of differentiable $(p,q)$-forms of $X$
with values in $E$ and $A_X^{p,q}(E)=\Gamma(X,\sA_X^{p,q}(E))$
the space of its global sections.

\item $T_X$ the holomorphic tangent bundle of $X$.

\item $K_{X}=A_X^{0,*}(T_X)=\oplus_{i} \Gamma(X,\sA^{0,i}_{X}(T_{X}))$ the
Kodaira-Spencer differential graded Lie algebra of $X$.

\item  $\Def(X)$ the deformation space of $X$,
i.e. the base space of the semiuniversal deformation of $X$.

\item $\mathbf{Art}$ the category of local artinian $\mathbb{C}$-algebras.

\item If $L$ is a  differential graded Lie algebra (DGLA), we denote by
\[\Def_L=\frac{\text{Maurer-Cartan}}{\text{Gauge action}}
\colon \mathbf{Art}\to \mathbf{Set}\]%
the associated deformation functor
(see \cite{GoMil2}, \cite{ManettiDGLA}, \cite{defomanifolds}
for  precise definition and properties).
\end{itemize}

In this paper we shall need several times 
the following result (for a proof see e.g. Theorem 5.71 of
\cite{defomanifolds}).

\begin{theorem}[Schlessinger-Stasheff \cite{SchSta}]\label{thm.fundamental}
Let $L\to M$ be a morphism of differential graded Lie algebras. Assume that:
\begin{enumerate}
\item  $H^0(L)\to H^0(M)$ is surjective.

\item $H^1(L)\to H^1(M)$ is bijective.

\item  $H^2(L)\to H^2(M)$ is injective.
\end{enumerate}
Then $\Def_L\to \Def_M$ is an isomorphism.\end{theorem}

If $X,Y$ are complex manifolds and $E\to X$, $F\to X$ vector
bundles, we denote  $E\boxtimes F=p^*E\otimes q^*F$, where
$p\colon X\times Y\to X$ and $q\colon X\times Y\to Y$ are the
projections. By K\"{u}nneth formula we have
\[ H^i(X\times Y, E\boxtimes F)=\mathop{\oplus}_j H^j(X,E)\otimes
H^{i-j}(Y,F).\]

\section{Twisted deformations of complex manifolds}

Let's recall the construction of the germ $\Def(X)$ for a compact
complex manifold $X$.\\
The starting point is the Kodaira-Spencer differential graded Lie algebra
\[K_{X}=A^{0,*}_{X}(T_{X}).\]
Fix a hermitian metric on $X$ and denote by
$d^*\colon A^{0,1}_{X}(T_{X})\to A^{0,0}_{X}(T_{X})$
the formal adjoint of $d$.
Then we consider suitable Sobolev completion (we may use the
appendix of \cite{DK} as instruction booklet)
\[ A^{0,i}_{X}(T_{X})\subset L^i\]
such that the  induced
maps of Hilbert spaces
\[ d\colon L^1\to L^{2},\qquad d^*\colon L^1\to L^{0},\qquad [~-,-]\colon L^1\otimes L^1\to
L^{2}.\]%
are defined and bounded.\\
Then (see e.g. \cite{GoMil2})  $\Def(X)$ is isomorphic to the germ at
0 of the analytic subvariety of $L^1$ defined by the equations
\[ dx+\frac{1}{2}[x,x]=0,\qquad d^*x=0.\]%
According to elliptic regularity, Hodge theory on compact manifolds and the results of
\cite{GoMil2}, the germ $\Def(X)$ is finite dimensional and well defined; in particular
different choices of  metric and Sobolev norms give isomorphic germs.

\begin{remark}
If the Kodaira-Spencer algebra $K_X$ is formal,
i.e. quasiisomorphic to its cohomology, then by \cite{GoMil1}, \cite{GoMil2},
the space $\Def(X)$ is analytically isomorphic to the germ at $0$ of the
nullcone of the quadratic map
\[H^1(X,T_X)\to H^2(X,T_X),\qquad \theta\mapsto [\theta,\theta].\]
In general the Kodaira-Spencer DGLA is not formal. In fact, it is a consequence of
\cite{Douady} that, if $X$ is the Iwasawa manifold  then $K_{X\times \mathbb{P}^1}$ is not formal;
more generally Ravi Vakil proved \cite{vakil}, putting together
the results of \cite{Mnev} and \cite{Diffeo}, that for every analytic
singularity $(U,0)$ defined  over $\mathbb{Z}$ there exists a complex surface $S$ with
very ample canonical bundle such that
$\Def(S)\cong (U\times \mathbb{C}^n,0)$ for some integer $n\ge 0$.\\
Choosing $U=\{(x,y)\in \mathbb{C}^2\mid xy(x-y)=0\}$ and taking $S$ as above,
the DGLA $K_S$ cannot be formal.
\end{remark}

Assume now that $(B,\wedge)$ is a finite dimensional graded
algebra of non negative degrees, say $B=\oplus_{i\ge 0}B^i$ and
$\dim_{\mathbb{C}}B<+\infty$.\\
Taking the natural extensions of the  operators $d,d^*$ and $[~,~]$ on $L\otimes B$,
\[ d_B(x\otimes h)=dx\otimes b,\quad d^*_B(x\otimes b)=
d^*x\otimes b,\quad [x\otimes b,y\otimes c]_B=
(-1)^{\bar{b}\;\bar{y}}[x,y]\otimes (b\wedge c),\]%
we may define an analytic germ
\[ \Def(X,B)=\{ x\in (L^1\otimes B^0)\oplus (L^0\otimes B^1)\mid
d_B^* x=0,\; d_Bx+\frac{1}{2}[x,x]_B=0\}.\]%

According to \cite[Thm's 3.9, 3.11]{GoMil2} (see also \cite[Thm. 4.7]{ManettiDGLA}),
the germ $\Def(X,B)$ is a hull (in the sense of \cite{Sch}) of the functor
\[ \Def_{K_X\otimes B}\colon \mathbf{Art}\to \mathbf{Set}\]
and then its  Zariski tangent space is isomorphic to
$(H^1(T_X)\otimes B^0)\oplus (H^0(T_X)\otimes B^1)$, while its
obstruction space is contained in
\[ H^2(K_X\otimes B)=
(H^2(T_X)\otimes B^0)\oplus (H^1(T_X)\otimes B^1)\oplus (H^0(T_X)\otimes B^2).\]
Note moreover that if $K_X$ is formal, then also $K_X\otimes B$ is formal.\\

Assume now that $B$ has a unit $1$ and $\dim B^0=1$; if $b_1,\ldots,
b_q$ is a basis of $B^1$, then $\Def(X,B)$ is the set of pairs
$(x,\sum_i y_i\otimes b_i)\in L^1\oplus (L^0\otimes B^1)$ such
that
\[ d^*x=0,\quad dx+\frac{1}{2}[x,x]=0,\quad
\sum_i (dy_i+[x,y_i])\otimes b_i=0,\quad \sum_{i<j}
[y_i,y_j]\otimes b_i\wedge b_j=0.\]%

The inclusion $\mathbb{C}=B^0\subset B$ induces a closed embedding
$\Def(X)\subset \Def(X,B)$, while the projection $B\to B^0$
induces an analytic retraction $r\colon \Def(X,B)\to \Def(X)$.\\
The fiber $r^{-1}(0)$ is the germ at $0$ of
\[ \left\{\sum_i y_i\otimes b_i\in L^0\otimes B^1\mid dy_i=0,\quad
\sum_{i<j}[y_i,y_j]\otimes b_i\wedge b_j=0\right\}.\]%
Since the kernel of $d\colon L^0\to L^1$ is $H^0(T_X)$, the above
set is equal to
\[ \{\sum_i y_i\otimes b_i\in H^0(T_X)\otimes B^1\mid
\sum_{i<j}[y_i,y_j]\otimes b_i\wedge b_j=0\}.\]%
In particular, if the product $\wedge\colon\bigwedge^2B^1\to B^2$
is injective, then $r^{-1}(0)$ is isomorphic to the commuting
variety
\[ C(q, H^0(T_X))=\{ (y_1,\ldots,y_q)\in H^0(T_X)^{\oplus q}\mid
[y_i,y_j]=0\text{ for every }i,j\}.\]

\bigskip
\section{Basic facts about commuting varieties}

Let $L$ be a finite dimensional complex Lie algebra and $q$ a
positive integer.\\
The affine scheme
\[ C(q,L)=\{(a_1,\ldots,a_q)\in L^{\oplus q}\mid [a_i,a_j]=0\;\text{ for every } i,j\}.\]
is called  the  $q$-th commuting
variety of $L$.
Clearly $C(q,L\oplus M)=C(q,L)\times C(q,M)$; while if $p\le q$,
then the projection on the first factors $C(q,L)\to C(p,L)$ is
surjective; in particular if $C(q,L)$ is irreducible, then also
$C(p,L)$ is irreducible.\\

The structure of the varieties $C(q,L)$ has been studied by
several people. The case  $L=\mathfrak{sl}(n,\mathbb{C})=H^0(T_{\mathbb{P}^{n-1}})$ has
been studied in Gerstenhaber \cite{Ger}; he proved  in
particular that $C(2,\mathfrak{sl}(n,\mathbb{C}))$ is irreducible  for every $n$
(this fact was also proved independently by Motzkin and Taussky
\cite{Motz}). It is a well known open (and hard) problem to determine
whether $C(2,\mathfrak{sl}(n,\mathbb{C}))$ (defined by the ideal generated by 
brackets) is a reduced scheme. 
Moreover, according to Richardson \cite{Rich}, the variety  $C(2,L)$ is irreducible for
every reductive Lie algebra $L$.

\begin{lemma}\label{lem.irred}
If the commuting variety $C(q,\mathfrak{sl}(n,\mathbb{C}))$ is irreducible
then 
\[ q< 3+\frac{8n-12}{(n-2)^2},\qquad \text{ for $n$ even},\]
\[ q< 3+\frac{8}{n-3},\qquad \text{ for $n$ odd}.\]
\end{lemma}

\begin{proof}
This proof  is based on the ideas of \cite{iarro}. Assume $n\ge
4$, $C(q,\mathfrak{sl}(n))$ irreducible and consider the projection on
the first factor
\[ \pi\colon C(q,\mathfrak{sl}(n,\mathbb{C}))\to
C(1,\mathfrak{sl}(n,\mathbb{C}))=\mathfrak{sl}(n,\mathbb{C}).\]%
Let $D\in \mathfrak{sl}(n,\mathbb{C})$ be a diagonal
matrix with distinct eigenvalues, then every matrix commuting with $D$ must be
diagonal. Therefore the fiber $\pi^{-1}(D)$ is irreducible of
dimension $(n-1)(q-1)$ and the dimension of $C(q,\mathfrak{sl}(n,\mathbb{C}))$ is
less than or equal
to $n^{2}-1+(n-1)(q-1)$.\\
On the other hand,  let $r$ be the integral part of $n/2$ and let
$N\subset \mathfrak{sl}(n,\mathbb{C})$ be the closed subset of matrices $A$ such
that $A^2=0$. It is easy to see that $N$ is irreducible of
dimension $2r(n-r)$ and therefore, for the generic $A\in N$, we
have
\[\dim\pi^{-1}(A)<  n^{2}-1-2r(n-r)+(n-1)(q-1).\]
After a possible change of basis, every $A\in N$ belongs to
the space
\[ H=\{(h_{ij})\in \mathfrak{sl}(n,\mathbb{C})\mid h_{ij}\not=0 \text{ only if }
i> r,\; j\le r \}.\]
Then $H$ is an abelian subalgebra of
$\mathfrak{sl}(n,\mathbb{C})$ and therefore $\{A\}\times H^{\oplus q-1}\subset
\pi^{-1}(A)$. In particular
\[ r(n-r)(q-1)=\dim H^{\oplus q-1}<  n^{2}-1+(n-1)(q-1)-2r(n-r).\]
A straightforward computation gives the inequalities of the lemma.
\end{proof}

\begin{remark}
If $\mathfrak{R}(\mathbb{Z}^q,\operatorname{PGL}(n,\mathbb{C}))$ is the space
of representations $\rho\colon \mathbb{Z}^q\to \operatorname{PGL}(n,\mathbb{C})$, then
the exponential gives an isomorphism between the germ at 0 of the commuting variety
$C(q,\mathfrak{sl}(n,\mathbb{C}))$ and a neigbourhood of the trivial representation in
$\mathfrak{R}(\mathbb{Z}^q,\operatorname{PGL}(n,\mathbb{C}))$.
\end{remark}

\bigskip
\section{Deformations of products}

Let $X,Y$ be two compact complex manifolds and denote by
\[ \sX\to \Def(X),\qquad \sY\to \Def(Y)\]
their semiuniversal deformations. The product $\sX\times \sY\to
\Def(X)\times \Def(Y)$ is a deformation of $X\times Y$ and
therefore induces a morphism of analytic germs
\[ \alpha\colon \Def(X)\times \Def(Y)\to \Def(X\times Y).\]

\begin{lemma}
The morphism $\alpha$ is a closed embedding;  it is an isomorphism
if and only if $H^0(T_X)\otimes H^1(\Oh_Y)=H^1(\Oh_X)\otimes
H^0(T_Y)=0$.\end{lemma}

\begin{proof}
By K\"{u}nneth formula
\[ H^1(T_{X\times Y})\simeq
H^1(T_X)\oplus H^1(T_Y)\oplus (H^0(T_X)\otimes H^1(\Oh_Y))\oplus
(H^1(\Oh_X)\otimes H^0(T_Y)).\]%
A morphism of analytic germs is a
closed embedding if and only if it is injective on Zariski tangent
spaces and the differential of $\alpha$ is equal to the natural
embedding
\[ H^1(T_X)\oplus H^1(T_Y)\to H^1(T_{X\times Y}).\]
The obstruction map associated to $\alpha$ is
\[ H^2(T_X)\oplus H^2(T_Y)\to H^2(T_{X\times Y})\]
and, again by K\"{u}nneth  formula,  it is injective. Therefore, if
$H^0(T_X)\otimes H^1(\Oh_Y)=H^1(\Oh_X)\otimes H^0(T_Y)=0$, then
the differential of $\alpha$ is bijective and $\alpha$ is an
isomorphism.
\end{proof}

The condition $H^0(T_X)\otimes H^1(\Oh_Y)=H^1(\Oh_X)\otimes
H^0(T_Y)=0$ is satisfied in most cases; for instance, a theorem
of Matsumura \cite{Matsu2} implies that $H^0(T_X)=0$ for
every compact manifold of general type $X$.\\
If $H^1(\Oh_X)\otimes H^0(T_Y)\not=0$,
then it is easy to describe deformations of $X\times Y$
that are not a product. Assume for simplicity $X$ K\"{a}hler, then
$b_1(X)\not=0$ and there exists at least one surjective homomorphism
$\pi_1(X)\mapor{g}\mathbb{Z}$. On the other hand, since $H^0(T_Y)\not=0$, there exists at least
a nontrivial one parameter subgroup
$\{\theta_t\}\subset\Aut(Y)$, $t\in \mathbb{C}$, of holomorphic automorphisms of $Y$.\\
Therefore we get a family of representations
\[ \rho_t\colon \pi_1(X)\to \Aut(Y),\qquad \rho_t(\gamma)=\theta_t^{g(\gamma)},\qquad t\in \mathbb{C}\]
inducing a family of locally trivial analytic $Y$-bundles over $X$.\\
Moreover,  Kodaira and Spencer  \cite{KS2}  proved that
projective spaces $\mathbb{P}^n$ and complex tori
$(\mathbb{C}^q/\Gamma)$ have unobstructed deformations, while
the product
$(\mathbb{C}^q/\Gamma)\times \mathbb{P}^n$
has obstructed deformations for every $q\ge 2$ and every $n\ge 1$ \cite[page 436]{KS2}.
This was the first example of obstructed manifold; the same example is
discussed, with a different approach, also in \cite{Douady}.

\begin{lemma}\label{lem.semiform}
For every pair of compact  K\"{a}hler manifolds
$X,Y$ there exists an injective quasi-isomorphism of differential
graded Lie algebras
\[ (B_{X}\otimes K_{Y})\oplus (B_{Y}\otimes K_{X})\to K_{X\times Y}.\]
\end{lemma}

\begin{proof}
Assume  $X,Y$  compact complex  manifolds and denote by
\[ p\colon X\times Y\to X,\qquad q\colon X\times Y\to Y\]
the projections. Since $T_{X\times Y}$ is the direct sum of
$p^{*}T_{X}=T_X\boxtimes\Oh_Y$ and $q^{*}T_{Y}=\Oh_X\boxtimes T_Y$
we have
\[ H^{i}(X\times Y,T_{X\times Y})=
H^{i}(X\times Y,p^{*}T_{X})\oplus  H^{i}(X\times Y,q^{*}T_{Y})\]
and, by K\"{u}nneth formula
\[ H^{i}(X\times Y,p^{*}T_{X})=
\mathop{\oplus}_{j}H^{j}(T_{X})\otimes H^{i-j}(\Oh_{Y})=\mathop{\oplus}_{j}H^{j}(T_{X})\otimes B_Y^{i-j},\]
\[ H^{i}(X\times Y,q^{*}T_{Y})=
\mathop{\oplus}_{j}H^{j}(\Oh_{X})\otimes H^{i-j}(T_{Y})=\mathop{\oplus}_{j}B_X^{j}\otimes H^{i-j}(T_{Y}).\]
The isomorphism $T_{X\times
Y}=p^{*}T_{X}\oplus q^{*}T_{Y}$ allows to define two injective
morphisms of differential graded Lie algebras
\[ p^{*}\colon K_{X}\to K_{X\times Y},\qquad
q^{*}\colon K_{Y}\to K_{X\times Y}.\]
We note that $p^{*}$ is injective in cohomology and the image of
$H^{i}(T_{X})$ is the subspace $H^{i}(T_{X})\otimes
H^{0}(\Oh_{Y})\subset H^{i}(X\times Y,p^{*}T_{X})$; similarly for
the morphism $q^{*}$.\\
Denote by $\bar{\Omega}^*_{X}$ the graded algebra of antiholomorphic differential
forms on $X$, endowed with the wedge product. More precisely
$\bar{\Omega}^*_{X}=\oplus \bar{\Omega}_{X}^{i}$ where
\[ \bar{\Omega}_{X}^{i}=\ker(\de\colon\Gamma(X,\sA^{0,i}_{X})\to
\Gamma(X,\sA^{1,i}_{X})).\]
If $X$ is compact K\"{a}hler, then  the
natural maps
\[\bar{\Omega}^*_{X}\cap \ker\bar{\de}\to \bar{\Omega}^*_{X},\qquad
\bar{\Omega}^*_{X}\cap \ker\bar{\de}\to
B_{X}=\frac{\ker\bar{\de}}{\Image\bar{\de}}\] are both
isomorphisms (see e.g. \cite{DGMS}) and therefore there exists an
isomorphism of graded algebras  $\bar{\Omega}^*_{X}\cong B_{X}$
independent from the K\"{a}hler metric. Denote also by
$\bar{\Omega}^*_{Y}$ the graded algebra of antiholomorphic
differential
forms on $Y$.\\
We can define two morphisms
\[ h\colon \bar{\Omega}^*_{Y}\otimes K_{X}\to K_{X\times Y}, \qquad
h(\phi\otimes \eta)=q^{*}(\phi)\wedge p^{*}(\eta).\]
\[ k\colon \bar{\Omega}^*_{X}\otimes K_{Y}\to K_{X\times Y}, \qquad
k(\psi\otimes \mu)=p^{*}(\psi)\wedge q^{*}(\mu).\]
It is straightforward to check that $h,k$ are morphisms of
differential graded Lie algebras and that the image of $h$ commutes
with the image of $k$. This implies that the map
\[ h\oplus k\colon (\bar{\Omega}^*_{X}\otimes K_{Y})\oplus (\bar{\Omega}^*_{Y}\otimes K_{X})
\to K_{X\times Y}\]
is a morphism of differential graded Lie algebras.
If $X$ and $Y$ are both K\"ahler manifolds, then according to
K\"{u}nneth formula the morphism $h\oplus k$ is a
quasiisomorphism.
\end{proof}

\begin{theorem}\label{thm.prodotto}
Let $X,Y$ be compact K\"{a}hler manifolds. Then
there exists an isomorphism of germs
\[ \Def(X\times Y)\simeq \Def(X,B_Y)\times\Def(Y,B_X).\]
Moreover if $K_X$ and $K_Y$ are formal differential graded Lie algebras, then also
$K_{X\times Y}$ is formal.
\end{theorem}

\begin{proof} By Artin's theorem on solutions of analytic
equations \cite{Artin68}, it is sufficient to show that there exists a formal isomorphism.\\
The germ $\Def(X\times Y)$ is a hull of the functor
$\Def_{K_{X\times Y}}$, while  the germ $\Def(X,B_Y)\times\Def(Y,B_X)$ is a hull of
the functor $\Def_{K_{X}\otimes B_{Y}}\times \Def_{K_{Y}\otimes B_X}$.\\
By Lemma \ref{lem.semiform} there exists a quasiisomorphism
$(B_{X}\otimes K_{Y})\oplus (B_{Y}\otimes K_{X})\to K_{X\times Y}$ 
and we can apply
Theorem \ref{thm.fundamental}.
\end{proof}

\begin{corollary}\label{cor.douady}
Let $\mathbb{C}^q/\Gamma$ be a complex torus of dimension $q$ and let $Y$ be a
compact K\"{a}hler manifold such that $H^1(\Oh_Y)=H^1(T_Y)=0$.
Then
\[\Def(\mathbb{C}^q/\Gamma\times Y)\simeq \mathbb{C}^{q^2}\times C(q,H^0(T_Y)).\]
Moreover, if $q>1$ then the following three conditions are
equivalent:
\begin{enumerate}

\item $\Def(\mathbb{C}^q/\Gamma\times Y)$ is smooth.

\item There exists the universal deformations of
$\mathbb{C}^q/\Gamma\times Y$.

\item The Lie algebra $H^0(Y,T_Y)$ is abelian.
\end{enumerate}
\end{corollary}

\begin{proof} The first part and the equivalence
\biimplica{1}{3} are immediate consequences of
Theorem \ref{thm.prodotto}.  We only prove the equivalence of conditions 2 and 3.\\
Let $z_1,\ldots, z_q$ be linear coordinates on $\mathbb{C}^q$. Denoting by
$\Cqbar=\langle d\bar{z_1},\ldots,d\bar{z_q}\rangle$ the space of
invariant $(0,1)$-form on the torus, there exists an isomorphism
$\bigwedge^*\Cqbar=B_{\mathbb{C}^q/\Gamma}$. Similarly the inclusion
\[
\bigwedge^*\Cqbar\otimes H^0(T_{\mathbb{C}^q/\Gamma})\hookrightarrow
K_{\mathbb{C}^q/\Gamma}\]%
is a quasiisomorphism of differential graded Lie algebras and then the inclusion
\[ Q:=\bigwedge^*\Cqbar\otimes \left(H^0(T_{\mathbb{C}^q/\Gamma})\oplus H^0(T_Y)\right)\hookrightarrow
K_{\mathbb{C}^q/\Gamma\times Y}\] is bijective in $H^0,H^1$ and injective in $H^2$.\\
In particular the functor $\Def_Q$ is isomorphic to the functors
of deformations of $\mathbb{C}^q/\Gamma\times Y$. Since $Q$ has trivial
differential, the functor $\Def_Q$ is prorepresentable if and only
if the gauge action is trivial or, equivalently, if and only if 
$[Q^0,Q^1]=0$. It is straightforward to check
that $[Q^0,Q^1]=0$ if and only if $H^0(T_Y)$
is abelian.
\end{proof}

\bigskip

\section{Deformations of $\mathbb{C}^q/\Gamma\times \mathbb{P}^n$}

Let $\mathbb{C}^q/\Gamma$ be a complex torus of dimension $q$ and $z_1,\ldots, z_q$ linear coordinates
on $\mathbb{C}^q$.
For later use and notational purposes, we  reprove the first part of
Corollary \ref{cor.douady} when $Y$ is the projective space $\mathbb{P}^n$.
In this case, according to K\"{u}nneth formula, the
inclusion
\begin{equation}\label{equ.Q}
Q:=\bigwedge^*\Cqbar\otimes (H^0(T_{\mathbb{C}^q/\Gamma})\oplus H^0(T_{\mathbb{P}^n}))
\hookrightarrow
K_{\mathbb{C}^q/\Gamma\times {\mathbb{P}^n}}
\end{equation}
is a quasiisomorphism of differential graded Lie algebras.

\begin{corollary}\label{cor.toripro}
The deformation space $\Def(\mathbb{C}^q/\Gamma\times \mathbb{P}^n)$ is analytically isomorphic to the
germ at 0 of
$\mathbb{C}^{q^2}\times C(q,\mathfrak{sl}(n+1))$. In particular
it is singular for every $n\ge 1$, $q\ge 2$ and it is
reducible for every $n\ge 3$ and $q\ge 3+8/(n-2)$.\end{corollary}

\begin{proof} Observe that
$H^0(T_{\mathbb{P}^n})\simeq \mathfrak{sl}(n+1)$ and apply the previous results.
\end{proof}

Our next goal is to prove a similar result
for deformations of polarized varieties. More precisely we assume that $\mathbb{C}^q/\Gamma$ is
an abelian variety, we fix an ample line bundle $L$ on $\mathbb{C}^q/\Gamma\times \mathbb{P}^n$ and we
consider deformations of  the pair $(\mathbb{C}^q/\Gamma\times \mathbb{P}^n, L)$.\\

We recall
that an Appell-Humbert (A.-H.) data  on a complex torus $\mathbb{C}^q/\Gamma$
is a pair $(\alpha, H)$,
where $H$ is  a hermitian form  on $\mathbb{C}^n$ such that its imaginary
part $E$ is integral on $\Gamma\times \Gamma$ and
\[ \alpha\colon \Gamma\to \operatorname{U}(1),\qquad
\alpha(\gamma_1+\gamma_2)=
\alpha(\gamma_1)\alpha(\gamma_2)(-1)^{E(\gamma_1,\gamma_2)}.\]%
Denote by $L(\alpha,H)$ the line bundle on $\mathbb{C}^q/\Gamma$ with
factor of automorphy \cite[pag. 4]{kempf}
\[A_{\gamma}(z)=\alpha(\gamma)e^{\pi
(H(z,\gamma)+H(\gamma,\gamma)/2)},\qquad \gamma\in\Gamma,\;
z\in\mathbb{C}^q.\]
 It is well known that every line bundle on $\mathbb{C}^q/\Gamma$ is
 isomorphic to $L(\alpha,H)$ for a unique Appell-Humbert data
 $(\alpha,H)$; moreover the first Chern class of
 $L(\alpha,H)$ is equal to the invariant $(1,1)$ form corresponding to $E$
 \cite[Lemma 3.5]{kempf}.
In particular two line bundles $L(\alpha_1,H_1)$ ,
$L(\alpha_2,H_2)$ have the same Chern class if and only if
$H_1=H_2$.\\
The same proof of the Appell-Humbert theorem given in
\cite{kempf}, with minor and straightforward modifications, shows
that every line bundle on $\mathbb{C}^q/\Gamma\times \mathbb{P}^n$ is
isomorphic to
\[L(\alpha,H,d):=L(\alpha,H)\boxtimes \Oh(d)\]%
for some A.-H. data $(\alpha,H)$ and some integer $d$.\\

Denote for simplicity $X=\mathbb{C}^q/\Gamma\times \mathbb{P}^n$.
The Atiyah extension of a line bundle $L(\alpha,H,d)$ is the short
exact sequence
\[ 0\mapor{} \Oh_X\mapor{}\sD(L(\alpha,H,d))\mapor{\sigma}
T_X\mapor{}0,\]%
where $\sD(L(\alpha,H,d))$ is the sheaf of first order
differential operator on $L(\alpha,H,d)$ and $\sigma$ is the
principal symbol. The induced morphism in cohomology
\[ H^1(T_X)=H^1(T_{\mathbb{C}^q/\Gamma})\otimes H^0(\Oh_{\mathbb{P}^n})\oplus
H^1(\Oh_{\mathbb{C}^q/\Gamma})\otimes H^0(T_{\mathbb{P}^n})\mapor{\contr c_1}
H^2(\Oh_X)\] is equal to the contraction with the first Chern
class of the line bundle $L(\alpha,H,d)$ (\cite{atiyah}).

\begin{lemma}\label{lem.ciuno}
Let $c_1$ be the first Chern
class of the line bundle $L(\alpha,H,d)$ on $X=\mathbb{C}^q/\Gamma\times \mathbb{P}^n$ 
and assume that $\det(H)\not=0$. Then $\contr
c_1(H^1(\Oh_{\mathbb{C}^q/\Gamma})\otimes H^0(T_{\mathbb{P}^n}))=0$ and
\[H^1(T_{\mathbb{C}^q/\Gamma})\otimes H^0(\Oh_{\mathbb{P}^n})\mapor{\contr c_1}
H^2(\Oh_X)\]%
is surjective; in particular $H^2( \sD(L(\alpha,H,d)))\mapor{}
H^2(T_X)$ is injective.
\end{lemma}

\begin{proof} The first part is clear since
$\contr c_1(H^1(\Oh_{\mathbb{C}^q/\Gamma})\otimes H^0(T_{\mathbb{P}^n}))\subset
H^1(\Oh_{\mathbb{C}^q/\Gamma})\otimes H^1(\Oh_{\mathbb{P}^n})$.\\
The map
\[H^1(T_{\mathbb{C}^q/\Gamma})\otimes H^0(\Oh_{\mathbb{P}^n})\mapor{\contr c_1}
H^2(\Oh_X)=H^2(\Oh_{\mathbb{C}^q/\Gamma})\otimes H^0(\Oh_{\mathbb{P}^n})\]%
can be written as $\contr c_1= e\otimes Id$, where $Id$ is the
identity on $H^0(\Oh_{\mathbb{P}^n})$ and
\[e\colon H^1(T_{\mathbb{C}^q/\Gamma})\to H^2(\Oh_{\mathbb{C}^q/\Gamma})\]%
is the contraction
with the first Chern class of $L(\alpha,H)$.\\
The elements of $H^2(\Oh_{\mathbb{C}^q/\Gamma})$ are represented by
invariant $(0,2)$-forms; more precisely, if $z_1,\ldots, z_q$ are
linear coordinates on $\mathbb{C}^q$, then a basis of
$H^2(\Oh_{\mathbb{C}^q/\Gamma})$ is given by the  forms $d\bar{z}_i\wedge
d\bar{z}_j$, for $i<j$. Similarly a basis of $
H^1(T_{\mathbb{C}^q/\Gamma})$ is given by the invariant tensors $
d\bar{z}_i\otimes \desude{~}{z_j}$, for $i,j=1,\ldots,q$.\\
The first Chern class of $L(\alpha,H)$ is given by the invariant
form $\sum h_{rs}dz_r\wedge d\bar{z}_s$, ,with $(h_{rs})$ a scalar
multiple of $H$, and
\[ e\left(d\bar{z}_i\otimes \desude{~}{z_j}\right)=
d\bar{z}_i\otimes \desude{~}{z_j}\contr
\sum_{r,s} h_{rs}dz_r\wedge d\bar{z}_s= \sum_s
h_{js}d\bar{z}_i\wedge d\bar{z}_s.\]%
Therefore the surjectivity of the matrix $(h_{rs})$
 implies the surjectivity of $e$.
\end{proof}

\begin{theorem}\label{thm.defline}

 Let $L$ be a line bundle on
$X=\mathbb{C}^q/\Gamma\times \mathbb{P}^{n}$ and denote by
$\Def(X, L)$ the deformation space of the pair
$(X, L)$.\begin{enumerate}

\item If $L\in \Pic^0(X)$, then the natural morphism $\Def(X, L)\to \Def(X)$
is smooth of relative dimension $q$.

\item If $L$ is ample then there exists a smooth morphism of analytic germs
\[ \Def(X, L)\to C(q,\mathfrak{sl}(n+1,\mathbb{C})).\]
\end{enumerate}
\end{theorem}

\begin{proof} The first part is a general result which is true
for every K\"{a}hler manifold $X$.
In fact by Deligne's theorem \cite{deligne},
for every deformation $\sX\to Spec(A)$ of $X$, with $A$ a local Noetherian
ring, the group $H^1(\Oh_\sX)$ is a free $A$-module of rank $h^1(\Oh_X)$ and then there
are no obstructions to
extend a topologically trivial line bundles over any deformation of $X$.\\
Assume now that
$L$ is ample; then  $L=L(\alpha,H,d)$
for some $d>0$ and $H$ positive definite. Consider the Dolbeault
resolution of the Atiyah extension of $L$
\[ 0\mapor{} A^{0,*}_X\mapor{}A^{0,*}_X(\sD(L))\mapor{\sigma}
A^{0,*}_X(T_X)\mapor{}0.\]%
The differential graded Lie algebra $A^{0,*}_X(\sD(L))$ governs the
deformations of the pair $(X,L)$ (see \cite{ClemensPetri}).\\
Denoting by $P$ the fiber product of $\sigma$ and the injective quasiisomorphism 
$Q\to A^{0,*}_X(T_X)$ of Equation (\ref{equ.Q}), we have a commutative diagram
\[\begin{array}{ccccccccc}
0&\mapor{}&A^{0,*}_X&\mapor{}&P&\mapor{}&Q&\mapor{}&0\\
&&\Vert&&\mapver{\phi}&&\mapver{}&&\\
0&\mapor{}&A^{0,*}_X&\mapor{}&A^{0,*}_X(\sD(L))&\mapor{\sigma}&
A^{0,*}_X(T_X)&\mapor{}&0\end{array}\]
with exact rows and vertical quasiisomorphisms.
In particular, according to Theorem \ref{thm.fundamental}, 
the functor $\Def_P$ is isomorphic to the
functor of deformations of the pair $(X,L)$.\\
The subspace $I=\bigwedge^1\Cqbar\otimes
H^0(T_{\mathbb{C}^q/\Gamma})\subset Q^1$ is an  ideal of the DGLA $Q$;
denote by $R=Q/I$ the corresponding quotient DGLA. Since
$R^1=\bigwedge^1\Cqbar\otimes H^0(T_{\mathbb{P}^{n}})$, the commuting
variety $C(q,\mathfrak{sl}(n+1,\mathbb{C}))$ is exactly the deformation space of
the functor $\Def_R$. Therefore it is sufficient to prove that the
morphism $P\to R$ is surjective on $H^1$ and injective on $H^2$;
this follows easily  from Lemma \ref{lem.ciuno}.
\end{proof}

\bigskip

\section{Natural deformations of complete intersections}

Let $X$ be a smooth complex manifold of dimension $n$ and let $D_1,\ldots,D_m$
be smooth divisors with $0<m\le n-2$. Assume that $D_1,\ldots ,D_m$ intersect
transversally
on a smooth subvariety $S$ of codimension $m$.\\
The  natural deformations of $S$ are the deformations obtained by
deforming $X$ and the divisors $D_i$. More precisely, let
$\Def_{X;D_1,\ldots,D_m}\colon \mathbf{Art}\to\mathbf{Set}$ be the functor of
infinitesimal deformations of  the holomorphic map
$\stackrel{\circ}{\cup}D_i\to X$. Equivalently every element of
$\Def_{X;D_1,\ldots,D_m}(A)$ is the data of a deformation $\sX\to
Spec(A)$ of $X$ and deformations $\sD_i\subset \sX$ of the
smooth hypersurfaces $D_i$. Since $\cap_i \sD_i$ is a deformation
of $S$ over $Spec(A)$, it is well defined a natural transformation
\[ \ope{Nat}\colon \Def_{X;D_1,\ldots,D_m}\to \Def_S.\]
The (infinitesimal) natural deformations of $S$ are the ones lying
in the image of $\ope{Nat}$.\\

The aim of this section is to give a sufficient condition for the completeness of the
natural deformations of $S$.
Since $S$ is complete intersection,  the ideal sheaf
$\sI_S\subset\Oh_X$ admits the Koszul resolution
\[ 0\mapor{}\Oh_X(-\sum_{i=1}^m D_i)\mapor{}\cdots\mapor{}
\mathop{\oplus}_{i<j}\Oh_X(-D_i-D_j)\mapor{}\mathop{\oplus}_{i}\Oh_X(-D_i)
\mapor{}\sI_S\mapor{}0.\]
Denote by $T_X(-\log D)\subset T_X$ the subsheaf of vector fields
which are tangent to $D_i$ for every $i$, and by
$N_{D_i|X}\simeq \Oh_{D_i}(D_i)$ the normal bundle of $D_i$ in $X$.
There exists a short exact double sequence:
\[\begin{array}{ccccccccc}
&&0&&0&&0&&\\
&&\mapver{}&&\mapver{}&&\mapver{}&&\\
0&\mapor{}&\sL&\mapor{}&T_X\otimes \sI_S&\mapor{}&\mathop{\oplus_i} \sI_S N_{D_i|X}
&\mapor{}&0\\
&&\mapver{}&&\mapver{}&&\mapver{}&&\\
0&\mapor{}&T_X(-\log D)&\mapor{}&T_X&\mapor{}&\mathop{\oplus_i} N_{D_i|X}&\mapor{}&0\\
&&\mapver{\alpha}&&\mapver{}&&\mapver{}&&\\
0&\mapor{}&T_S&\mapor{}&T_X\otimes\Oh_S&\mapor{}&N_{S|X}&\mapor{}&0\\
&&\mapver{}&&\mapver{}&&\mapver{}&&\\
&&0&&0&&0&&\end{array}\]

\begin{theorem}\label{thm.costability}
For every $A\subset\{1,\ldots,m\}$ denote by $D_A=\sum_{i\in A}D_i$. Assume that
\begin{enumerate}

\item $H^{|A|+1}(T_X(-D_A))=0$ for every $A\not=\emptyset$.

\item $H^{|A|}(\Oh_X(D_i-D_A))=0$ for every $i=1,\ldots,m$ and every
$A\not=\{i\},\emptyset$.
\end{enumerate}
Then the natural transformation
$\ope{Nat}\colon \Def_{X;D_1,\ldots,D_m}\to \Def_S$ is smooth.
\end{theorem}

\begin{proof}
Considering the cohomology of the tensor product of $T_X$ with the Koszul resolution of
$\sI_S$, we get immediately that
the first condition implies that $H^2(T_X\otimes\sI_S)=0$.\\
Let $i$ be fixed; assume that $A\not=\emptyset$ and $i\not\in A$, then from the exact
sequence
\[ 0\mapor{}\Oh_X(D_i-D_{A\cup\{i\}})\mapor{}\Oh_X(D_i-D_{A})
\mapor{}\Oh_{D_i}(D_i-D_{A})\mapor{}0\]
we get $H^{|A|}(\Oh_{D_i}(D_i-D_A))=0$.
Considering the cohomology of the tensor product of $\Oh_{D_i}(D_i)$ with the Koszul
resolution of the ideal of $S$ in $D_i$ we get that
$H^1(\sI_S N_{D_i|X})=0$. Therefore $H^2(\sL)=0$ and the
morphism $\alpha\colon T_X(-\log D)\to T_S$ is surjective in $H^1$ and
injective in $H^2$.\\
By general results of deformation theory
(see e.g. \cite{namba})
tangent and obstruction spaces of the functor $\Def_{X;D_1,\ldots,D_m}$ are
$H^1(T_X(-\log D))$ and $H^2(T_X(-\log D))$ respectively. Therefore
$\ope{Nat}$ is surjective on tangent spaces, injective on obstruction spaces
and therefore it is smooth.
\end{proof}

\begin{remark}
For $m=1$ the Theorem \ref{thm.costability} 
reduces to a particular case of Horikawa's costability theorem.
\end{remark}

\bigskip
\section{Some new example of obstructed irregular surfaces}

In order to apply Theorem \ref{thm.costability} to the variety
$X=\mathbb{C}^q/\Gamma\times \mathbb{P}^n$ we need the determination of cohomology groups
of the line bundles $L(\alpha,H,d)$. Notice  that $K_X=L(0,0,-n-1)$
and
\[ L(\alpha_1,H_1,d_1)\otimes L(\alpha_2,H_2,d_2)=L(\alpha_1\alpha_2,H_1+H_2,d_1+d_2).\]

\begin{lemma}\label{lem.detcoh}
Let  $X=\mathbb{C}^q/\Gamma\times \mathbb{P}^n$, then
\begin{enumerate}
\item If $\alpha\not=1$ then $H^i(L(\alpha,0,d))=0$ for every
$i,d\in\mathbb{Z}$.

\item If $H$ is positive definite and $d\ge -n$, then
$H^i(L(\alpha,H,d))=0$ for every
$i>0$.

\item If $H$ is negative definite and $d\le -2$, then
$H^i(L(\alpha,H,d))=H^i(T_X\otimes L(\alpha,H,d))=0$ for every
$i\le q+n-2$.

\item If $H$ is negative definite and $d\le -n-2$, then
$H^{q+n-1}(T_X\otimes L(\alpha,H,d))=0$.
\end{enumerate}
\end{lemma}

\begin{proof} The determination of the cohomology of line bundles
on  $\mathbb{C}^q/\Gamma$ \cite[Th. 3.9]{kempf} gives that:
\begin{itemize}

\item  If $\alpha\not=1$ then $H^i(L(\alpha,0))=0$ for every $i,$.

\item If $H$ is negative definite, then $H^i(L(\alpha,H))=0$ for
every $i<q$.
\end{itemize}
By K\"{u}nneth formula we get  item (1). Assume now that $H$ is
negative definite; then by K\"{u}nneth formula we get
$H^i(L(\alpha,H,d))=0$ for every $\alpha$, every $d<0$ and every
$i<q+n$;  Serre duality implies (2).\\
The  bundle $T_X\otimes L(\alpha,H,d)$ is the direct sum of
$L(\alpha,H)\boxtimes T_{\mathbb{P}^n}(d)$ and $q$ copies of
$L(\alpha,H,d)$. Since $H^i(T_{\mathbb{P}^n}(d))=0$ for every $d\le
-2$ and every $i\le n-2$ we get item (3).\\
If $d\le -n-2$ then $H^{n-1}(T_{\mathbb{P}^n}(d))=0$ and this implies
item (4).
\end{proof}

\begin{definition}
Let $D$ be a divisor of a section of the line bundle
$L(\alpha,H,d)$. We shall call the pair $(H,d)$ the \emph{homology
type} of $D$.\end{definition}

The above definition is justified since the Poincar\'e dual of $D$
is the sum of $d$ times the hyperplane class on $\mathbb{P}^n$ and the
class represented by the imaginary part of $H$.

\begin{proposition}\label{prop.abecommu}
 Assume $q,n,d\ge 2$ and $q+n+d\ge 7$.
Let $S$ be a smooth ample hypersurface of $\mathbb{C}^q/\Gamma\times
\mathbb{P}^{n-1}$ of homology type $(H,d)$. 
Then every deformation of
$S$ is projective and $\Def(S)$ has the same singularity type of the
commuting variety $C(q,\mathfrak{sl}(n,\mathbb{C}))$.
\end{proposition}

\begin{proof}
The hypersurface $S$ is the zero locus of a section $s$ of a line
bundle $L=L(\alpha,H,d)$ for some semi-character $\alpha$. Since
$S$ is ample the hermitian form $H$ is positive definite and then,
according to Theorem \ref{thm.defline} there exists a smooth
morphism
\[ \Def(X,L)\to C(q,\mathfrak{sl}(n,\mathbb{C})).\]
Since $H^1(L)=0$, the section $s$ extends to every deformation of
the pair $(X,L)$ and then the natural morphism $\Def(X,S)\to
\Def(X,L)$ is smooth. The ampleness of $S$ also implies that every
deformation of the pair $(X,S)$ is projective.\\
On the other hand, by Lemma \ref{lem.detcoh}  $H^2(T_X(-S))=0$ and
then by Theorem \ref{thm.costability} the morphism $\Def(X,S)\to
\Def(S)$ is smooth.
\end{proof}

\begin{remark} The fact that obstructions to smoothness of
$\Def(X,S)\to \Def(X,L)$ are contained  in $H^1(L)$ is well known and easy to
prove directly. However, it is interesting to see this fact also
as a consequence of the exact sequence
\[ 0\to T_X(-\log S)\mapor{\alpha}\sD(L)\mapor{\beta}L\to 0,\]
\[ \alpha(\chi)=\chi-s^{-1}\chi(s),\qquad \beta(\phi)=\phi(s).\]
Notice that $\alpha$ is a morphism of sheaves of  Lie algebras.
\end{remark}

\begin{example}\label{exa.easyexample} 
Let $A$ be an abelian surface and
$S$ a smooth surface of general type contained in $A\times
\mathbb{P}^1$. Then $\Oh(S)=L(\alpha,H,d)$ for some $H$ positive definite and $d\ge 3$;
according to Proposition \ref{prop.abecommu}
$\Def(S)$ has the singularity type of $C(2,\mathfrak{sl}(2,\mathbb{C}))$.
Since two matrices in $\mathfrak{sl}(2,\mathbb{C})$ commute if and only
if they are linearly dependent, the commuting variety
$C(2,\mathfrak{sl}(2,\mathbb{C}))$ is equal to determinantal variety of matrices
$2\times 3$ of rank $\le 1$.\end{example}

\begin{theorem}\label{thm.main}
Let $S$ be a smooth
complete intersection of $m=q+n-3$ smooth divisors
$D_1,\ldots,D_m$ of $X=\mathbb{C}^q/\Gamma\times
\mathbb{P}^{n-1}$.
Assume that $\Oh(D_i)=L(\alpha_i,H,d)$ with
$H$ is positive definite,
$d\ge 2$,  $d(q+n-3)\ge n+1$
and $\alpha_i\not=\alpha_j$ for every $i\not=j$.\\
Then every deformation of $S$ is
projective and $\Def(S)$ has the same singularity type of the
commuting variety $C(q,\mathfrak{sl}(n,\mathbb{C}))$.\\
In particular, if  $q,n\ge 2$ then
$S$ is obstructed and $\Def(S)$  is
reducible for  $n\ge 4$ and $q\ge 3+8/(n-3)$.
\end{theorem}

\begin{proof}
By adjunction formula the surface $S$ has ample canonical bundle and then
every deformation of $S$ is projective.\\
By Lemma \ref{lem.detcoh}   and
Theorem \ref{thm.costability} the morphism
$\ope{Nat}\colon \Def(X;D_1,\ldots,D_m)\to \Def(S)$ is smooth.
According to Theorem \ref{thm.defline} there exists a smooth
morphism
\[ \Def(X,L(\alpha_1,H,d))\to C(q,\mathfrak{sl}(n,\mathbb{C})).\]
Therefore it is sufficient to prove that  that
the morphism
\[ \Def(X;D_1,\ldots,D_m)\to \Def(X;\Oh(D_1))\]
is smooth.
First of all, since $H^1(D_i)=0$ for every $i$ the morphism
\[ \Def(X;D_1,\ldots,D_m)\to \Def(X;\Oh(D_1),\ldots,\Oh(D_m))\]
is smooth. On the other side, since $\Oh(D_i-D_j)\in \Pic^0(X)$, by Theorem
\ref{thm.defline} also the morphism
\[ \Def(X;\Oh(D_1),\ldots,\Oh(D_m))\to \Def(X;\Oh(D_1))\]
is smooth.
\end{proof}

\bigskip

\end{document}